\documentclass[12pt,leqno]{amsart}
\usepackage[latin9]{inputenc}
\usepackage{verbatim}
\usepackage{amsthm}
\usepackage{amstext}
\usepackage{amssymb}
\usepackage[dvipdfmx]{graphicx}
\usepackage{xcolor}
\usepackage{soul}

\usepackage{mathrsfs}
\usepackage{amscd}

\makeatletter
\numberwithin{equation}{section}

\usepackage{amsthm}\usepackage{amscd}\usepackage{bm}

\usepackage[dvipdfmx,pdfstartview=FitH,pdfborderstyle={/S/B/W 1}]{hyperref}

\theoremstyle{plain}
\newtheorem{main theorem}{Main Theorem}
\newtheorem{theorem}{Theorem}[section]
\newtheorem{lemma}[theorem]{Lemma}
\newtheorem{conjecture}[theorem]{Conjecture}
\newtheorem{corollary}[theorem]{Corollary}
\newtheorem{proposition}[theorem]{Proposition}
\newtheorem{claim}[theorem]{Claim}

\theoremstyle{definition}

\newtheorem{remark}[theorem]{Remark}

\newtheorem{problem}[theorem]{Problem}

\newcommand{\diam}{\mathrm{diam}}
\newcommand{\supp}{\mathrm{supp}}

\newcommand{\ind}{\mathrm{ind}}
\newcommand{\coind}{\mathrm{coind}}
\newcommand{\conn}{\mathrm{conn}}
\makeatother

\begin{document}

\title[$G$-index, topological dynamics and marker property]{$G$-index, topological dynamics and marker property}

\author{Masaki Tsukamoto, Mitsunobu Tsutaya, Masahiko Yoshinaga}



\date{\today}

\thanks{M. Tsukamoto was supported by JSPS KAKENHI 18K03275.  M. Tsutaya was supported by 
JSPS KAKENHI 19K14535. M. Yoshinaga was supported by JSPS KAKENHI 18H01115.}

\maketitle

\begin{abstract}
Given an action of a finite group $G$, we can define its index. 
The $G$-index roughly measures a size of the given $G$-space.
We explore connections between the $G$-index theory and topological dynamics.
For a fixed-point free dynamical system, we study the $\mathbb{Z}_p$-index of the set of 
$p$-periodic points.
We find that its growth is at most linear in $p$.
As an application, we construct a free dynamical system which does not have the marker property.
This solves a problem which has been open for several years. 
\end{abstract}

\section{Introduction}   \label{section: introduction}

\subsection{Background on $\mathbb{Z}_p$-index}  \label{subsection: background on Z_p-index}

When a finite group $G$ acts on a topological space, we can define its \textit{index}, which roughly measures 
a ``size'' of the given $G$-space.
The $G$-index theory is a generalization of the Borsuk--Ulam theorem, and it has several striking applications to 
combinatorics \cite{Matousek}.

The purpose of this paper is to exhibit new applications of the $G$-index to \textit{dynamical systems theory}.
In particular, we prove the following theorem. 
The terminologies of dynamical systems will be explained in \S \ref{subsection: dynamical systems and Z_p index} and 
\S \ref{subsection: marker property}.

\begin{theorem}  \label{theorem: marker property}
There exists a free dynamical system which does not have the marker property.
\end{theorem} 

This solves a problem which has been open for several years.
We explain a background of this theorem in \S \ref{subsection: marker property}.

First we introduce basic definitions of the $G$-index theory. 
Our presentation follows 
a book of Matou$\mathrm{\check{s}}$ek \cite[Sections 6.1 and 6.2]{Matousek}.
In this paper, we consider only the case of $G = \mathbb{Z}/p\mathbb{Z}$ for prime numbers $p$.
We set $\mathbb{Z}_p := \mathbb{Z}/p\mathbb{Z}$.
A pair $(X, T)$ is called a \textbf{$\mathbb{Z}_p$-space} if $X$ is a topological space and 
$T:X\to X$ is a homeomorphism satisfying $T^p = \mathrm{id}$.
We often omit $T$ from the notation and simply say ``$X$ is a $\mathbb{Z}_p$-space''. 

Let $(X, T)$ be a $\mathbb{Z}_p$-space.
It is said to be \textbf{free} if $T^a x \neq x$ for any $1\leq a \leq p-1$ and any $x\in X$.
Since $p$ is a prime number, this is equivalent to the condition that $T x\neq x$ for any $x\in X$.
Let $n\geq 0$ be an integer.
A free $\mathbb{Z}_p$-space $(X, T)$ is called an \textbf{$E_n\mathbb{Z}_p$-space} if it satisfies 
\begin{itemize}
  \item $X$ is an $n$-dimensional finite simplicial complex.
  \item $T:X\to X$ is a simplicial map, i.e. mapping each simplex to a simplex affinely.
  \item $X$ is $(n-1)$-connected, i.e. the $k$-th homotopy group $\pi_k (X)$ vanishes for all $0\leq k \leq n-1$.
\end{itemize}

For example, the discrete space $\mathbb{Z}_p$ with the natural $\mathbb{Z}_p$ action is an $\mathbb{E}_0\mathbb{Z}_p$-space.
(We consider that $\mathbb{Z}_p$ is $(-1)$-connected.) In general, the join of $(n+1)$-copies of $\mathbb{Z}_p$ is an 
$\mathbb{E}_n\mathbb{Z}_p$-space:
\[  E_n\mathbb{Z}_p = \underbrace{\mathbb{Z}_p *\mathbb{Z}_p*\dots *\mathbb{Z}_p}_{\text{$(n+1)$-times}}.  \]
Here $\mathbb{Z}_p$ acts on each component of the join simultaneously. 
(If a finite simplicial complex $X$ is $n$-connected and a finite simplicial complex $Y$ is $m$-connected then 
the join $X*Y$ is $(n+m+2)$-connected; \cite[4.4.3 Proposition]{Matousek}. )

There also exist other models of $E_n\mathbb{Z}_p$-space.
For example,
\[  E_n\mathbb{Z}_p = \begin{cases}  
                                S^n & (\text{when $n$ is odd}) \\
                                S^{n-1}*\mathbb{Z}_p & (\text{when $n$ is even})
                               \end{cases}. \]
So $E_n\mathbb{Z}_p$-space is not unique. 
However they are ``essentially'' unique for our purpose because 
any $E_n\mathbb{Z}_p$-space equivariantly and continuously maps to any other (\cite[6.2.2 Lemma]{Matousek}).

Let $(X, T)$ be a free $\mathbb{Z}_p$-space.
We define 
\[ \ind_p (X, T) = \min\left\{n\geq 0\middle|\, \exists f:X\to E_n\mathbb{Z}_p: 
   \text{$\mathbb{Z}_p$-equivariant continuous mapping}\right\}.\] 
We set $\ind_p (X, T) = \infty$ if there is no equivariant continuous map $f:X\to E_n\mathbb{Z}_p$ for any $n\geq 0$.
We often abbreviate $\ind_p(X,T)$ as $\ind_p X$.
We have $\ind_p E_n\mathbb{Z}_p = n$ (\cite[6.2.5 Theorem]{Matousek}).
This can be seen as a generalization of the Borsuk--Ulam theorem.
We use the convention that if $X$ is empty then $\ind_p X = -1$.

For a (non-empty) topological space $X$, 
we define its \textbf{connectivity} $\conn(X)$ as the smallest integer $n\geq -1$ satisfying 
$\pi_{n+1}(X) \neq 0$.
We set $\conn(\emptyset) := -2$.
It is known (\cite[6.2.4 Proposition]{Matousek}) that for a $\mathbb{Z}_p$-space $(X,T)$
\begin{equation}  \label{eq: index and connectivity}
   \ind_p (X,T) \geq \conn (X) +1.
\end{equation}

\subsection{Dynamical systems and $\mathbb{Z}_p$-index}  \label{subsection: dynamical systems and Z_p index}

Next we introduce basic definitions of dynamical systems and state our first main result.
A pair $(X, T)$ is called a \textbf{dynamical system} if $X$ is a compact metrizable space and 
$T:X\to X$ is a homeomorphism\footnote{Notice that we assume the compactness of $X$.
This is essential for our results.}. 
We sometimes omit $T$ from the notation and simply say ``$X$ is a dynamical system'' when 
the action $T$ can be understood in the context.

Let $(X, T)$ be a dynamical system. For a natural number $n=1,2,3,\dots$, we define 
$P_n(X, T)$ as the set of $n$-periodic points:
\[ P_n(X, T) := \left\{x\in X\middle|\, T^n x = x\right\}. \]
We often abbreviate $P_n(X, T)$ as $P_n(X)$.

A dynamical system is said to be \textbf{free} if it has no periodic points, namely $P_n(X, T) = \emptyset$
for all $n\geq 1$.
A dynamical system is said to be \textbf{fixed-point free} if it has no fixed point, namely $P_1(X, T) = \emptyset$.

Suppose a dynamical system $(X, T)$ has no fixed-point. 
For each prime number $p$, we consider a $\mathbb{Z}_p$-space
\[  \left(P_p(X, T), T\right). \]
Since $(X, T)$ has no fixed point, this is a free $\mathbb{Z}_p$-space.
So we can consider its $\mathbb{Z}_p$-index, and
we get a sequence 
\[  \ind_p \left(P_p(X, T), T\right), \quad (p=2,3,5,7,11,13,17, \dots). \]
Now we ask a question: \textit{What kind of sequence can appear? Is there any restriction?}

Indeed it has a restriction. Its growth is at most linear in $p$:

\begin{theorem} \label{theorem: first main result}
Let $(X, T)$ be a fixed-point free dynamical system.
Then $\ind_p P_p(X)$ is finite for all prime numbers $p$, and the sequece 
\[  \ind_p P_p(X) \quad (p=2,3,5,7, 11, \dots) \]
has at most linear growth. 
Namely there exists a positive number $C$ satisfying 
\[ \ind_p P_p(X) < C\cdot p \]
for all prime numbers $p$.
\end{theorem}

We have 
\[  \conn P_p(X) \leq \ind_p P_p(X) -1 \]
by (\ref{eq: index and connectivity}).
So we get the next corollary.

\begin{corollary}
Let $(X, T)$ be a fixed-point free dynamical system.
Then the connectivity of $P_p(X)$ grows at most linearly in $p$.
\end{corollary}

\subsection{Marker property} \label{subsection: marker property}

Theorem \ref{theorem: first main result} indicates that we can detect a hidden structure of dynamical systems
by using $\mathbb{Z}_p$-index theory.
But is there any application of this new structure?
Surprisingly, it has an application to a problem which seems to have nothing to do with $\mathbb{Z}_p$-index theory.

A dynamical system $(X, T)$ is said to have the \textbf{marker property} if for any natural number $N$
there exists an open set $U\subset X$ satisfying the following two conditions.
\begin{itemize} 
   \item No point in $U$ returns to $U$ within $N$ steps, namely 
   \[   U \cap T^{-n} U = \emptyset  \quad \text{for all $1\leq n\leq N$}. \] 
   \item Every orbit of $T$ has non-empty intersection with $U$, namely
\[   X = \bigcup_{n\in \mathbb{Z}} T^{-n} U.   \]
\end{itemize}
For example, infinite minimal dynamical systems\footnote{A dynamical system $(X,T)$ is called an \textbf{infinite minimal system} 
if $X$ is an infinite set and every orbit of $T$ is dense in $X$.
For example, an irrational rotation of the circle is an infinite minimal system.}
and their extensions have the marker property.

It is easy to see that if a dynamical system $(X, T)$ has the marker property then it is free. 
Indeed, suppose it has a periodic point $p$ of period $N\geq 1$.
Take an open set $U\subset X$ satisfying $U\cap T^{-n}U = \emptyset$ for all $1\leq n\leq N$ and 
$X = \bigcup_{n\in \mathbb{Z}} T^{-n}U$.
Then $q:=T^n(p) \in U$ for some $n$.
The point $q$ is also an $N$-periodic point. But $q = T^N(q) \not \in U$. 
This is a contradiction.

The marker propery has been intensively used in the context of \textit{mean dimension theory}
\cite{Lindenstrauss, Gutman Jaworski, Gutman periodic points, Gutman--Lindenstrauss--Tsukamoto, Gutman--Qiao--Tsukamoto, 
Lindenstrauss--Tsukamoto double VP, Gutman--Tsukamoto, Tsukamoto mdim potential}.
Mean dimension is a topological invariant of dynamical systems which counts the number of parameters of 
dynamical systems per iterate.
It was introduced by Gromov \cite{Gromov}. See also the paper of Lindenstrauss--Weiss \cite{Lindenstrauss--Weiss}.
We will not use mean dimension in the main body of the paper. So we do not provide the precise definition.

Here are two samples of theorems using the marker property in their assumptions
(\cite[Proposition 6.14]{Lindenstrauss}\footnote{\cite[Proposition 6.14]{Lindenstrauss} was stated only for 
extensions of infinite minimal systems. But actually the proof is valid for all dynamical systems with the marker property.} 
and \cite[Main Theorem 1]{Gutman--Qiao--Tsukamoto}):

\begin{theorem}[\cite{Lindenstrauss}]  \label{theorem: Lindenstrauss 1999}
Suppose a dynamical system $(X, T)$ has the marker property.
Then $(X, T)$ has zero mean dimension if and only if it is isomorphic to a projective limit of 
finite entropy dynamical systems.
\end{theorem}

\begin{theorem}[\cite{Gutman--Qiao--Tsukamoto}] \label{theorem: GQT 2019}
Suppose a dynamical system $(X, T)$ has the marker property. Let $N$ be a natural number.
If $(X, T)$ has mean dimension smaller than $N/2$ then it can be embedded in the shift action on 
the infinite dimensional cube
\[  \left([0,1]^N\right)^{\mathbb{Z}} = \cdots \times [0,1]^N\times [0,1]^N\times [0,1]^N \times \cdots. \]
\end{theorem}

As we remarked in the above, if a dynamical system has the marker property then it is free.
It is very natural to ask whether the converse holds or not:

\begin{problem} \label{problem: marker property}
Does every free dynamical system have the marker property?
\end{problem}

Besides its naturalness, this question is also important for better understanding the range of applicability of 
various theorems using the marker property, such as the above
Theorems \ref{theorem: Lindenstrauss 1999} and \ref{theorem: GQT 2019}.

Problem \ref{problem: marker property} was stated by Gutman 
(\cite[Problem 5.4]{Gutman Jaworski} and \cite[Problem 3.4]{Gutman periodic points}).
The origin of the question goes back to the paper of Lindenstrauss \cite{Lindenstrauss}.
In that paper he introduced a certain topological analogue of the Rokhlin tower lemma
(\cite[Lemma 3.3]{Lindenstrauss})\footnote{We will review this in Lemma \ref{lemma: Lindenstrauss} below.}.
This Rokhlin-type lemma is an equivalent form of the marker property. 
After proving the lemma, he wrote \cite[p. 234, lines 1 and 2]{Lindenstrauss}:
\begin{quote}
 In this form, it does not seem that this Rokhlin-type Lemma can be extended to a more general setup.
\end{quote}
So he more or less conjectured a negative answer to Problem \ref{problem: marker property}.
Probably many serious readers of \cite{Lindenstrauss} had an interest in the question.

Problem \ref{problem: marker property} looks an innocent question.
However, unexpectedly, it turns out to be difficult.
The main reason is that \textit{it has an infinite dimensional nature}.
Gutman \cite[Theorem 6.1]{Gutman Jaworski} proved that every \textit{finite dimensional} free
dynamical system has the marker property.
Here a dynamical system $(X, T)$ is said to be \textbf{finite dimensional} if the topological dimension 
(Lebesgue covering dimension) of $X$ is finite.
So the main objects in the study of Problem \ref{problem: marker property} are infinite dimensional dynamical systems.
The problem is difficult simply because our current technology for understanding infinite dimensional systems 
is very limited.
Probably the most important discovery of the paper is that $\mathbb{Z}_p$-index theory provides a new 
tool for studying infinite dimensional dynamical systems.

As we already mentioned in \S \ref{subsection: background on Z_p-index}, 
the answer to Problem \ref{problem: marker property} is \textit{no}.
This is our second main result.

\begin{theorem}[$=$ Theorem \ref{theorem: marker property}] \label{theorem: marker property revisted}
There exists a free dynamical system which does not have the marker property.
\end{theorem}

The relation between this theorem and Theorem \ref{theorem: first main result} in \S 
\ref{subsection: dynamical systems and Z_p index} is roughly as follows.
In the proof of Theorem \ref{theorem: first main result} we introduce a certain \textit{universal}
fixed-point free dynamical system $\mathcal{X}(N,\delta)$, which is \textit{larger} than the given 
$(X,T)$ from the viewpoint of $\mathbb{Z}_p$-index theory.
It turns out that this universal system also becomes a basic building block for constructing a 
free dynamical system which does not have the marker property.

\section{Proof of Theorem \ref{theorem: first main result}} \label{section: proof of first main result}

In this section we introduce dynamical systems $\mathcal{X}(N,\delta)$ which have a certain universal property.
Theorem \ref{theorem: first main result} immediately follows from the universality of $\mathcal{X}(N,\delta)$.
The universality will be also a crucial ingredient for the proof of Theorem \ref{theorem: marker property} later.

\subsection{$\varepsilon$-embedding}  \label{subsection: varepsilon-embedding}

This subsection is a preparation.
Let $(X, d)$ and $(Y, d')$ be metric spaces. We assume that $X$ is compact.
Let $\varepsilon>0$. A continuous map $f:X\to Y$ is called an \textbf{$\varepsilon$-embedding}
if for all $y\in Y$
\[ \diam f^{-1}(y) < \varepsilon. \]
Since we assume the compactness of $X$, if $f$ is an $\varepsilon$-embedding then there exists $\delta>0$
such that for any two points $x_1, x_2\in X$
\[  d'\left(f(x_1), f(x_2)\right) < \delta \> \Longrightarrow  \>  d(x_1, x_2) < \varepsilon. \]

\begin{lemma} \label{lemma: varepsilon-embedding}
  Let $\varepsilon>0$ and let $(X, d)$ be a compact metric space.
There are a natural number $N$ and an $\varepsilon$-embedding from $X$ to the $N$-dimensional cube $[0,1]^N$.
\end{lemma}

\begin{proof}
We can assume $\diam X \leq 1$ (by using some scale change).
We take open balls $B_1, \dots, B_N$ of radius $\varepsilon/2$ covering $X$.
Let $x_1, \dots, x_N$ be the centers of $B_1, \dots, B_N$ respectively.
We define 
\[  f: X\to [0,1]^N, \quad x\mapsto \left(d(x, x_1), d(x, x_2), \dots, d(x,x_N)\right). \]
This is an $\varepsilon$-embedding. Indeed, suppose $x, y\in X$ satisfy $f(x) = f(y)$.
Suppose, say, $x\in B_1$. Then $d(y, x_1) = d(x, x_1) < \varepsilon/2$ and hence $y\in B_1$.
Since $B_1$ is an $\varepsilon/2$-ball, we get $d(x, y) < \varepsilon$.
\end{proof}

\subsection{Universal fixed-point free dynamical systems} \label{subsection: universal fixed-point free dynamical systems}

Let $N$ be a natural number. We consider the infinite product of the copies of the $N$-dimensional cube:
\[ \left([0,1]^N\right)^{\mathbb{Z}} = \cdots [0,1]^N\times [0,1]^N\times [0,1]^N\times \cdots. \]
This becomes a dynamical system under the shift map 
\[ \sigma: \left([0,1]^N\right)^{\mathbb{Z}}\to \left([0,1]^N\right)^{\mathbb{Z}}, \quad 
     \left(x_n\right)_{n\in \mathbb{Z}} \mapsto \left(x_{n+1}\right)_{n\in \mathbb{Z}}. \]
Let $\delta>0$. We define a subsystem $\mathcal{X}(N, \delta)$ of $\left([0,1]^N\right)^{\mathbb{Z}}$ by 
\[ \mathcal{X}(N,\delta) := \left\{(x_n)_{n\in \mathbb{Z}}\middle|\, 
    \forall n\in \mathbb{Z}: |x_n-x_{n+1}| \geq \delta\right\}. \]
The pair $\left(\mathcal{X}(N, \delta), \sigma\right)$ becomes a dynamical system.
(We often omit $\sigma$ from the notation and simply write ``the dynamical system $\mathcal{X}(N, \delta)$''.)
This system is the central object of the paper.
It immediately follows from the definition that $\mathcal{X}(N,\delta)$ has no fixed point:
\[ P_1 \left(\mathcal{X}(N,\delta),\sigma\right) = \emptyset. \]

$\mathcal{X}(N,\delta)$ has the following \textit{universality}.

\begin{theorem}[\textbf{Universality of $\mathcal{X}(N,\delta)$}] \label{theorem: universality}
Let $(X, T)$ be a fixed-point free dynamical system.
There are a natural number $N$ and a positive number $\delta$ such that
there exists an equivariant continuous map from $X$ to 
$\mathcal{X}(N,\delta)$.
\end{theorem}

\begin{proof}
Let $d$ be a metric on $X$ compatible with the topology.
Since $(X,T)$ has no fixed point, there is a positive number $\varepsilon$ satisfying 
\[ \forall x\in X:  \> d\left(x, T x\right) > \varepsilon. \]
By Lemma \ref{lemma: varepsilon-embedding}, there is an $\varepsilon$-embedding 
$f:X\to [0,1]^N$ for some natural number $N$.
Since $X$ is compact, we can find $\delta>0$ such that for any two points $x, y\in X$
\[  \left|f(x)-f(y)\right| < \delta \> \Longrightarrow \> d(x, y) < \varepsilon. \]
Then for all $x\in X$
\begin{equation}  \label{eq: f(x) and f(Tx) are not close}
   \left|f(x)-f(Tx)\right| \geq \delta.
\end{equation}  

We define $F:X\to \left([0,1]^N\right)^{\mathbb{Z}}$ by 
\[  F(x) = \left(f(T^n x)\right)_{n\in \mathbb{Z}}. \]
$F$ is continuous and equivariant, i.e. $F\circ T = \sigma \circ F$.
By the above (\ref{eq: f(x) and f(Tx) are not close}), $F(x) \in \mathcal{X}(N,\delta)$ for all $x\in X$.
Thus $F$ is an equivariant continuous map from $X$ to $\mathcal{X}(N,\delta)$.
\end{proof}

Although the above proof is very simple, Theorem \ref{theorem: universality} is a key result of the paper.

\subsection{$P_p \left(\mathcal{X}(N,\delta)\right)$ and the proof of Theorem \ref{theorem: first main result}}
\label{subsection: proof of first main result}

Let $N$ be a natural number and $\delta$ a positive number.
For a prime number $p$, we study the set of $p$-periodic points of the dynamical system $\mathcal{X}(N,\delta)$.
From the definition, $P_p\left(\mathcal{X}(N,\delta)\right)$ is naturally identified with
\begin{equation}  \label{eq: P_p X(N,delta)}
   \left\{\left(x_0, x_1, \dots, x_{p-1}\right)\in \left([0,1]^N\right)^{\mathbb{Z}_p}\middle|\, 
       \forall n\in \mathbb{Z}_p: |x_n-x_{n+1}|\geq \delta\right\}.
\end{equation}
(Here the condition $|x_n-x_{n+1}|\geq \delta$ becomes $|x_{p-1}-x_0|\geq \delta$ when $n=p-1$.)
The group $\mathbb{Z}_p$ acts on this space by the cyclic shift:
\[  (x_0, x_1,x_2, \dots, x_{p-1}) \mapsto (x_1, x_2, \dots, x_{p-1}, x_0). \]
$P_p\left(\mathcal{X}(N,\delta)\right)$ is a free $\mathbb{Z}_p$-space.
From (\ref{eq: P_p X(N,delta)}), 
\[  P_p\left(\mathcal{X}(N,\delta)\right) \subset 
     \left(\mathbb{R}^N\right)^{\mathbb{Z}_p} \setminus \left\{(x, x, \dots, x)\middle|\, x\in \mathbb{R}^N\right\}. \]
The right-hand side is $\mathbb{Z}_p$-equivariantly homotopic to the $(Np-N-1)$-dimensional sphere 
\[  \left\{(x_0, x_1, \dots, x_{p-1})\in \left(\mathbb{R}^N\right)^{\mathbb{Z}_p}\middle|\, 
      x_0+x_1+\dots+x_{p-1}=0, \> |(x_0, x_1, \dots, x_{p-1})| = 1\right\}. \]
This is an $E_{Np-N-1}\mathbb{Z}_p$-space.
So its $\mathbb{Z}_p$-index is equal to $Np-N-1$.
Therefore 
\[ \ind_p P_p\left(\mathcal{X}(N,\delta)\right) \leq Np-N-1. \]
Now we are ready to prove Theorem \ref{theorem: first main result}.
For the convenience of readers, we write the statement again:

\begin{theorem}[$=$ Theorem \ref{theorem: first main result}]
Let $(X, T)$ be a fixed-point free dynamical system.
Then $\ind_p P_p(X)$ has at most linear growth in $p$.
\end{theorem}

\begin{proof}
From Theorem \ref{theorem: universality}, there is an equivariant continuous map 
$f:X\to \mathcal{X}(N,\delta)$ for some $N\geq 1$ and $\delta>0$.
Let $p$ be any prime number. Restricting $f$ to the set of $p$-periodic points, we get 
a $\mathbb{Z}_p$-equivariant continuous map 
\[ f:P_p(X) \to P_p\left(\mathcal{X}(N,\delta)\right). \]
Then 
\[ \ind_p P_p(X) \leq \ind_p P_p\left(\mathcal{X}(N,\delta)\right)  \leq Np-N-1. \]
\end{proof}

\section{Coindex of free $\mathbb{Z}_p$-spaces}  \label{section: coindex of free Z_p spaces}

The next three sections are preparations for constructing a free dynamical system which does not 
have the marker property.
In this section we introduce \textit{coindex} of free $\mathbb{Z}_p$-spaces and study its basic properties.
Most of the results in this section are certainly well-known (for example, see \cite[p.99]{Matousek}).
But we provide full proofs for the completeness.

Throughout this section we assume that $p$ is a prime number.
Let $(X, T)$ be a free $\mathbb{Z}_p$-space.
We define its \textbf{coindex} by 
\[ \coind_p (X, T) := \max\left\{n\geq 0\middle|\, \exists f: E_n\mathbb{Z}_p\to X: \text{ equivariant continuous map}\right\}. \]
We often abbreviate $\coind_p(X,T)$ as $\coind_p X$.
We use the convention that the coindex of the empty set is $-1$.

It is known that if there exists an equivariant continuous map from $E_m\mathbb{Z}_p$ to $E_n\mathbb{Z}_p$ then
$m\leq n$ (\cite[6.2.5 Theorem]{Matousek}). So we always have $\coind_p (X, T) \leq \ind_p(X,T)$.
We have $\coind_p E_n \mathbb{Z}_p = n$.

The reasons why we use $\coind_p X$ instead of $\ind_p X$ are the properties (2) and (3) of the next proposition.

\begin{proposition}[\textbf{Basic properties of coindex}] \label{prop: basic properties of coindex}
Let $X$ and $Y$ be free $\mathbb{Z}_p$-spaces.
   \begin{enumerate}
      \item If there exists an equivariant continuous map $f:X\to Y$ then $\coind_p X \leq \coind_p Y$.
      \item The product $X\times Y$ also becomes a free $\mathbb{Z}_p$-space. ($\mathbb{Z}_p$ acts on each component 
simultaneously.) Its coindex is given by 
     \[  \coind_p (X\times Y) = \min\left(\coind_p X , \coind_p Y\right). \] 
       \item The join\footnote{The join $X*Y$ is defined by 
    $X*Y = [0,1]\times X\times Y/\sim$, where the equivalence relation is given by 
\[  (0, x, y) \sim (0, x, y'), \quad (1, x, y) \sim (1, x', y) \]
for any $x, x'\in X$ and $y, y'\in Y$. The equivalence class of $(t, x, y)$ is usually written as $(1-t)x\oplus ty$.
The group action on the join $X*Y$ is given by the simultaneous actions on the component $X$ and the component $Y$.
We use the convention that if $Y=\emptyset$ then $X*Y = X$.} 
$X*Y$ is also a free $\mathbb{Z}_p$-space,
and we have 
\[  \coind_p (X*Y) \geq \coind_p X + \coind_p Y + 1. \] 
   \end{enumerate}
\end{proposition}

\begin{proof}
(1) Suppose $\coind_p X\geq n$. Then there exists an equivariant continuous map $g:E_n\mathbb{Z}_p\to X$.
The composition $f\circ g: E_n\mathbb{Z}_p\to Y$ is also equivariant and continuous.
So $\coind_p Y \geq n$.

(2) The projections from $X\times Y$ to each factors $X$ and $Y$ are equivariant continuous maps.
So, from (1), we get 
\[  \coind_p (X\times Y) \leq \min\left(\coind_p X, \coind_p Y \right). \]
Suppose we are given equivariant continuous maps 
\[  f:E_m \mathbb{Z}_p \to X, \quad g:E_n\mathbb{Z}_p\to Y. \]
Assume $m\leq n$. Then there exists an equivariant continuous map 
$h:E_m\mathbb{Z}_p\to E_n\mathbb{Z}_p$.
(This is obvious if we use the models 
$E_m \mathbb{Z}_p = \left(\mathbb{Z}_p\right)^{*(m+1)}\subset E_n\mathbb{Z}_p = \left(\mathbb{Z}_p\right)^{*(n+1)}$.)
Then 
\[ E_m \mathbb{Z}_p\to X\times Y, \quad u\mapsto \left(f(u), g(h(u))\right) \]
is an equivariant continuous map. So $m\leq \coind_p (X\times Y)$. Namely 
\[  \coind_p(X\times Y) \geq  \min\left(\coind_p X, \coind_p Y \right). \]

(3) Suppose we are given equivariant continuous maps
\[  f:E_m\mathbb{Z}_p\to X, \quad g:E_n\mathbb{Z}_p\to Y. \]
Then the join of the maps
\[ f*g:E_m\mathbb{Z}_p*E_n\mathbb{Z}_p \to X*Y, \quad 
    (1-t)u\oplus t v   \mapsto (1-t)f(u)\oplus t g(v)  \]
is also an equivariant continuous map.
If we use the model $E_n\mathbb{Z}_p = \left(\mathbb{Z}_p\right)^{*(n+1)}$ then 
\[  E_m\mathbb{Z}_p *E_n\mathbb{Z}_p = \left(\mathbb{Z}_p\right)^{*(m+n+2)} = E_{m+n+1}\mathbb{Z}_p. \]
Therefore $f*g$ is an equivariant continuous map from $E_{m+n+1}\mathbb{Z}_p$ to $X*Y$. So 
\[ \coind_p(X*Y) \geq m+n+1. \]
Namely 
\[ \coind_p (X*Y) \geq \coind_p X + \coind_p Y + 1. \] 
\end{proof}

Recall that a free $\mathbb{Z}_p$-space $(X, T)$ is called an $E_n\mathbb{Z}_p$-space if 
$X$ is an $n$-dimensional and $(n-1)$-connected finite simplicial complex with a simplicial map 
$T:X\to X$.
It follows from this defintion that if $(X, T)$ is an $E_n\mathbb{Z}_p$-space then, for any natural number $1\leq a\leq p-1$, 
the $\mathbb{Z}_p$-space $(X, T^a)$ is also an $E_n\mathbb{Z}_p$-space.
We use this fact in the proof of the next lemma.

\begin{lemma}  \label{lemma: coindex and iteration}
Let $(X, T)$ be a free $\mathbb{Z}_p$-space. For any natural number $1\leq a\leq p-1$
\[  \coind_p \left(X, T^a\right) = \coind_p (X, T). \]
\end{lemma}

\begin{proof}
Let $\left(E_n\mathbb{Z}_p, S\right)$ be an $E_n\mathbb{Z}_p$-space $(n\geq 0)$, and suppose we are given 
an equivariant continuous map $f:\left(E_n\mathbb{Z}_p, S\right) \to (X, T)$.
Then the same map $f$ also gives an equivariant continuous map from $\left(E_n\mathbb{Z}_p, S^a\right)$
to $\left(X, T^a\right)$. Since $\left(E_n\mathbb{Z}_p, S^a\right)$ is also an $E_n\mathbb{Z}_p$-space, 
we get $\coind_p\left(X, T^a\right) \geq n$. Namely 
\[  \coind_p\left(X, T^a\right) \geq \coind_p(X, T). \]
Take a natural number $1\leq b \leq p-1$ with $ab \equiv 1 \, (\mathrm{mod}.\, p)$.
Applying the same argument to $\left(X, T^a\right)$, we have 
\[  \coind_p\left(X, \left(T^a\right)^b\right) \geq \coind_p\left(X, T^a\right). \]
Since $T^{ab} = T$, the left-hand side is equal to $\coind_p(X, T)$.
Thus 
\[ \coind_p \left(X, T\right) \geq  \coind_p \left(X, T^a\right). \]
\end{proof}

Let $N$ be a natural number. We consider the product of $p$ copies of the $N$-dimensional cube:
\[  \left([0,1]^N\right)^{\mathbb{Z}_p} 
  = \underbrace{[0,1]^N\times [0,1]^N\times\cdots \times [0,1]^N}_{\text{$p$ times}}. \]
This becomes a (non-free) $\mathbb{Z}_p$-space under the cyclic shift:
\[ \sigma: \left([0,1]^N\right)^{\mathbb{Z}_p}  \to \left([0,1]^N\right)^{\mathbb{Z}_p}, \quad 
   \left(x_0, x_1, x_2, \dots, x_{p-1}\right) \mapsto \left(x_1, x_2, \dots, x_{p-1}, x_0\right). \]

\begin{corollary}  \label{cor: coindex of two important spaces}
Let $1\leq m\leq p-1$ be a natural number and let $\delta$ be a positive number.
We consider the following two $\mathbb{Z}_p$-subspaces of $\left([0,1]^N\right)^{\mathbb{Z}_p}$:
\begin{equation*}
  \begin{split}
   X &:= \left\{(x_0, x_1, \dots, x_{p-1})\in \left([0,1]^N\right)^{\mathbb{Z}_p}\middle|\, 
   \forall n\in \mathbb{Z}_p: |x_n-x_{n+1}|\geq \delta\right\}, \\
   Y &:= \left\{(x_0, x_1, \dots, x_{p-1})\in \left([0,1]^N\right)^{\mathbb{Z}_p}\middle|\,
    \forall n\in \mathbb{Z}_p: |x_n-x_{n+m}|\geq \delta\right\}.
  \end{split}
\end{equation*}
Then $(X,\sigma)$ and $(Y,\sigma)$ are free $\mathbb{Z}_p$-spaces, and 
\[ \coind_p (X, \sigma) = \coind_p (Y, \sigma). \]
\end{corollary}

\begin{proof}
We can immediately see from the definition that $(X, \sigma)$ and $(Y,\sigma)$ are free.
Take a natural number $1\leq l\leq p-1$ with 
$lm \equiv 1\, (\mathrm{mod}. \, p)$.
We define continuous maps $f:X\to Y$ and $g:Y\to X$ by 
\begin{equation*}
   \begin{split}
    f\left(x_0, x_1, \dots, x_{p-1}\right) = \left(x_0, x_l, x_{2l}, \dots, x_{l(p-1)}\right), \\
    g\left(y_0, y_1, \dots, y_{p-1}\right) = \left(y_0, y_m, y_{2m}, \dots, y_{m(p-1)}\right).   
   \end{split}
\end{equation*}
More precisely, 
\[ f\left((x_n)_{n\in \mathbb{Z}_p}\right) = \left(x_{ln}\right)_{n\in \mathbb{Z}_p}, \quad 
    g\left((y_n)_{n\in \mathbb{Z}_p}\right) =  \left(y_{mn}\right)_{n\in \mathbb{Z}_p}. \]
Then 
\[  f\circ g = \mathrm{id}, \quad g\circ f = \mathrm{id}. \]
Moreover 
\[  f\circ \sigma = \sigma^m \circ f, \quad g\circ \sigma^m = \sigma \circ g. \]
Therefore the $\mathbb{Z}_p$-space $(X, \sigma)$ is isomorphic to $(Y, \sigma^m)$.
Then 
\[  \coind_p (X,\sigma) = \coind_p\left(Y,\sigma^m\right) = \coind_p(Y,\sigma). \]
The last equality follows from Lemma \ref{lemma: coindex and iteration}.
\end{proof}

\section{Properties of $\mathcal{X}(N,\delta)$ and its variants} \label{section: X(N,delta) and its variants}

In this section we study a propety of $\mathcal{X}(N,\delta)$ introduced in 
\S \ref{subsection: universal fixed-point free dynamical systems} again.
We also introduce its variants.

\subsection{Join of dynamical systems}  \label{subsection: join of dynamical systems}
 
Let $(X, T)$ and $(Y,S)$ be dynamical systems.
We consider their join:
\[  (X*Y, T*S). \]
As usual, this is defined by 
\[ X*Y = [0,1]\times X\times Y/\sim, \quad \text{where }  (0, x, y)\sim (0, x, y'), \> (1, x, y) \sim (1, x', y), \]
\[ T*S\left((1-t)x\oplus t y\right) = (1-t) T(x)\oplus tS(y). \]
Here $(1-t)x\oplus ty$ denotes the equivalence class of $(t, x, y)$.
The join $(X*Y, T*S)$ is also a dynamical system. For any natural number $n$
\[  P_n(X*Y) = P_n(X)*P_n(Y). \]
In particular if $(X, T)$ and $(Y, S)$ are both fixed-point free, then their join is also fixed-point free.
Let $p$ be a prime number. Applying Proposition \ref{prop: basic properties of coindex} (3) to 
\[  P_p(X*Y) = P_p(X)*P_p(Y), \]
we have 
\begin{equation}  \label{eq: coind of join of dynamical system}
  \coind_p P_p(X*Y) \geq \coind_p P_p(X) + \coind_p P_p(Y) + 1.
\end{equation}
We are going to use this when $Y$ is a \textit{symbolic subshift}.

Conisder the \textbf{full-shift} on the alphabet $\{1,2,3\}$:
\[  \{1,2,3\}^\mathbb{Z} = \cdots\times \{1,2,3\}\times \{1,2,3\}\times \{1,2,3\}\times \cdots, \]
\[ \sigma: \{1,2,3\}^\mathbb{Z} \to \{1,2,3\}^\mathbb{Z}, \quad 
\sigma\left((x_n)_{n\in \mathbb{Z}}\right) = (x_{n+1})_{n\in \mathbb{Z}}. \]
The pair $\left(\{1,2,3\}^\mathbb{Z},\sigma\right)$ is a dynamical system.
We define its subsystem $\Sigma$ by 
\[\Sigma = \left\{(x_n)_{n\in \mathbb{Z}}  \in \{1,2,3\}^\mathbb{Z} \, \middle|\, 
                  \forall n\in \mathbb{Z}: x_n\neq x_{n+1}\right\}. \]

\begin{lemma}  \label{lemma: symbolic subshift}
  $(\Sigma, \sigma)$ has no fixed point. For any $m\geq 2$, $P_m(\Sigma) \neq \emptyset$.
\end{lemma}

\begin{proof}
$P_1(\Sigma) = \emptyset$ is obvious. Let $m\geq 2$. 
We are going to prove that $\Sigma$ has $m$-periodic points.
When $m$ is an even number, the point 
\[  \dots 12121212\dots \]
is an $m$-periodic point. 
Let $m=2l+1$ be odd with $l\geq 1$.
Define a word $u$ of length $m$ by 
\[  u = \underbrace{1212\dots 12}_{\text{$l$ times}} 3. \]
Then the point 
\[ \dots uuu\dots \]
is an $m$-periodic point.
\end{proof}

\begin{corollary}  \label{cor: increasing coindex by one}
 Let $(X, T)$ be a fixed-point free dynamical system.
There are a natural number $N$ and a positive number $\delta$ such that for all prime numbers $p$
\[  \coind_p P_p\left(\mathcal{X}(N,\delta)\right) \geq \coind_p P_p(X) + 1. \] 
\end{corollary}

\begin{proof}
Let $\left(\Sigma,\sigma\right)$ be the symbolic subshift introduced in the above. We consider the join 
\[  \left(X*\Sigma, T*\sigma\right). \]
This is fixed-point free. By (\ref{eq: coind of join of dynamical system})
\[  \coind_p P_p\left(X*\Sigma\right) \geq \coind_p P_p(X) + \coind_p P_p\left(\Sigma\right) + 1.\]
By Lemma \ref{lemma: symbolic subshift}, $P_p\left(\Sigma\right)$ is a non-empty finite set with a free 
$\mathbb{Z}_p$-action.
So it is an $E_0 \mathbb{Z}_p$-space. 
(Recall that the discrete space $\mathbb{Z}_p$ is an $E_n\mathbb{Z}_p$-space.)
Hence 
\[ \coind_p P_p\left(\Sigma\right) = 0. \]
Therefore 
\[   \coind_p P_p\left(X*\Sigma\right) \geq \coind_p P_p(X) + 1. \] 

Now we apply Theorem \ref{theorem: universality} (the universality of $\mathcal{X}(N,\delta)$)
to the join $X*\Sigma$.
There are $N\geq 1$ and $\delta>0$ such that there exists an equivariant continuous map 
\[  f: X*\Sigma  \to \mathcal{X}(N,\delta). \]
Let $p$ be any prime number.
Restricting $f$ to the set of $p$-periodic points, we have a $\mathbb{Z}_p$-equivariant continuous map
\[  f:   P_p\left(X*\Sigma\right)  \to P_p\left(\mathcal{X}(N,\delta)\right). \]
By Proposition \ref{prop: basic properties of coindex} (1)
\[ \coind_p P_p\left(\mathcal{X}(N,\delta)\right) 
    \geq \coind_p P_p\left(X*\Sigma\right) \geq \coind_p P_p(X) + 1. \] 
\end{proof}

\subsection{Variants of $\mathcal{X}(N,\delta)$}  \label{subsection: variants of X(N,delta)}

Let $N$ be a natural number.
We consider 
\[  \left([0,1]^N\right)^{\mathbb{Z}} = \cdots \times [0,1]^N\times [0,1]^N\times [0,1]^N \times \cdots. \]
This becomes a dynamical system under the shift 
\[ \sigma\left((x_n)_{n\in \mathbb{Z}}\right) = (x_{n+1})_{n\in \mathbb{Z}}. \]
Let $\delta$ be a positive number and $m$ a natural number.
We define a subsystem of $\left(\left([0,1]^N\right)^{\mathbb{Z}}, \sigma\right)$ by 
\[ \mathcal{X}_m(N,\delta) := \left\{(x_n)_{n\in \mathbb{Z}}\middle|\, \forall n\in \mathbb{Z}:
    |x_n-x_{n+m}| \geq \delta\right\}. \]
This is a fixed-point free dynamical system.
When $m=1$, it coincides with the system $\mathcal{X}(N,\delta)$:
\[ \mathcal{X}_1(N,\delta) = \mathcal{X}(N,\delta). \]

\begin{lemma}  \label{lemma: coind of X_m(N,delta)}
The dynamical system $\left(\mathcal{X}_m(N,\delta),\sigma\right)$ has no $m$-periodic point.
For any prime number $p>m$, 
\[  \coind_p P_p\left(\mathcal{X}_m(N,\delta)\right) = \coind_p P_p\left(\mathcal{X}(N,\delta)\right). \]
\end{lemma}

\begin{proof}
Obviously $P_m\left(\mathcal{X}_m(N,\delta)\right) = \emptyset$.
The set $P_p\left(\mathcal{X}_m(N,\delta)\right)$ is identified with 
\[ \left\{(x_0, x_1, \dots, x_{p-1})\in \left([0,1]^N\right)^{\mathbb{Z}_p}\middle|\, 
    \forall n\in \mathbb{Z}_p: |x_n-x_{n+m}| \geq \delta\right\}. \]
From Corollary \ref{cor: coindex of two important spaces}, for prime numbers $p>m$
\[  \coind_p P_p\left(\mathcal{X}_m(N,\delta)\right) = \coind_p P_p\left(\mathcal{X}(N,\delta)\right). \]
\end{proof}

\section{Lindenstrauss' lemma and its consequence}  \label{section: Lindentrauss' lemma and its consequence}

In this section we review a topological analogue of the Rokhlin tower lemma introduced by 
Lindenstrauss \cite{Lindenstrauss}.
We also study its consequence.

\subsection{Lindestrauss' lemma}  \label{subsection: Lindenstrauss' lemma}

As we mentioned in \S \ref{subsection: marker property}, 
Lindenstrauss \cite[Lemma 3.3]{Lindenstrauss} introduced a topological dynamics version of the Rokhlin tower lemma in ergodic theory.
Here is that lemma:

\begin{lemma}[\cite{Lindenstrauss}]  \label{lemma: Lindenstrauss}
Let $(X, T)$ be a dynamical system having the marker property. For any natural number $N$ there is 
a continuous function $\varphi:X \to \mathbb{R}$ such that the set 
\[  E := \{x\in X|\, \varphi(Tx) \neq \varphi(x)+1\} \]
satisfies $E\cap T^{-n} E = \emptyset$ for all $1\leq n\leq N$.
\end{lemma}

\begin{proof}
We explain Lindenstrauss' ingenious proof for the completeness\footnote{\cite[Lemma 3.3]{Lindenstrauss} was stated 
only for extensions of infinite minimal systems. 
But the proof is valid for all dynamical systems having the marker property.}.
From the definition of the marker property in \S \ref{subsection: marker property}, there is an open set $U\subset X$ satisfying 
\[  U\cap T^{-n}U = \emptyset \quad (\forall 1\leq n\leq N), \quad  X = \bigcup_{n\in \mathbb{Z}} T^n U. \]
Since $X$ is compact, there is a natural number $M\geq N$ satisfying 
\[  X = \bigcup_{n=0}^M T^n U. \]
Then we can find a compact set $K\subset U$ satisfying\footnote{This can be seen as follows.
Let $d$ be a metric on $X$. Set $U_n = T^n U$.
Since $X = \bigcup_{n=0}^M U_n$, we have $\max_{0\leq n\leq M} d(x, U_n^c) >0$ for all $x\in X$.
Since $X$ is compact, there is $\delta>0$ satisfying $\max_{0\leq n\leq M} d(x, U_n^c) \geq \delta$ for all $x\in X$.
We define a compact set $K_n\subset U_n$ by $K_n = \{x|\, d(x, U_n^c) \geq \delta\}$. Then 
$X = K_0\cup K_1\cup\dots \cup K_M$. Set $K := K_0\cup T^{-1}K_1\cup T^{-2}K_2\cup \dots \cup T^{-M}K_M$.
This $K$ satisfies the requirement.}
\[  X = \bigcup_{n=0}^M T^n K. \]
Take a continuous function $w:X\to [0,1]$ satisfying 
\[  w(x) = 1 \quad (\forall x\in K), \quad \supp\, w\subset U. \]

We consider a \textbf{Markovian random walk} on $X$ defined by 
\begin{itemize}
  \item When a particle is at a point $x$, the random walk ends with probability $w(x)$ and it moves to $T^{-1} x$
with probability $\left(1-w(x)\right)$.
\end{itemize}
Since $X = \bigcup_{n=0}^M T^n K$ and $w=1$ on $K$, this random walk stops or enters the set $\{w=1\}$ (and stops)
within $M$ steps.

For $x\in X$, we define $\varphi(x)$ as the \textbf{expected number of steps in the random walk starting at $x$}.
The explicit formula is as follows\footnote{Readers can easily understand this by checking the cases of some small $M$.
For example, if $M=1$ then 
\[ \varphi(x) = 1-w(x). \]
If $M=2$ then 
\[ \varphi(x) = \left(1-w(x)\right)w(T^{-1}x) + 2\left(1-w(x)\right) \left(1-w(T^{-1}x)\right) w(T^{-2}x). \]}:
\[ \varphi(x) = \sum_{n=1}^M n \cdot \left(\prod_{k=0}^{n-1} \left(1-w(T^{-k}x)\right)\right)\cdot w(T^{-n}x). \]   

Now we have defined a continuous function $\varphi:X\to \mathbb{R}$.
We are going to prove that the set $E= \{x|\, \varphi(Tx) \neq \varphi(x)+1\}$ satisfies the requirement.
Take a point $x\not \in U$. Then $w(x)=0$ and hence the random walk at $x$ moves to $T^{-1}x$ with probability one.
So $\varphi(x) = \varphi(T^{-1}x) + 1$. Therefore we have 
\[  E \subset T^{-1}U. \]
Then for $1\leq n\leq N$
\[ E\cap T^{-n} E \subset T^{-1}U \cap T^{-n-1}U = T^{-1}\left(U\cap T^{-n}U\right) = \emptyset. \]
\end{proof}

\begin{remark}
As is pointed out  by Gutman \cite[Theorem 7.3]{Gutman Jaworski}, 
the existence of a continuous function $\varphi$
stated in Lindenstrauss' lemma is indeed equivalent to the marker property.
Namely, \textit{a dynamical system $(X, T)$ has the marker property if and only if 
for any natural number $N$ there is a continuous function $\varphi:X \to \mathbb{R}$ such that 
the set $E = \{x\in X|\, \varphi(Tx)\neq \varphi(x)+1\}$ satisfies 
$E\cap T^{-n} E= \emptyset$ for all $1\leq n\leq N$.}
So Lindenstrauss' lemma provides an equivalent condition for the marker property.
\end{remark}

\subsection{Dynamical systems $\mathcal{Y}$ and $\mathcal{Z}$} 
\label{subsection: dynamical systems Y and Z}

In this subsection we introduce two dynamical systems related to the marker property.
We consider a circle $\mathbb{R}/2\mathbb{Z}$. (Here we use $\mathbb{R}/2\mathbb{Z}$ instead of
more natural $\mathbb{R}/\mathbb{Z}$ for later convenience.)
We define a metric $\rho$ on it by 
\[ \rho(x, y) := \min_{n\in \mathbb{Z}} |x-y-2n|. \]
The diameter of $\left(\mathbb{R}/2\mathbb{Z}, \rho\right)$ is one. 
(The distance between antipodal points is equal to one, e.g. $\rho(0,1) = 1$.)

We consider the infinite product of the copies of $\mathbb{R}/2\mathbb{Z}$:
\[   \left(\mathbb{R}/2\mathbb{Z}\right)^{\mathbb{Z}} = \cdots \times \mathbb{R}/2\mathbb{Z} \times \mathbb{R}/2\mathbb{Z} \times 
      \mathbb{R}/2\mathbb{Z} \times \cdots. \]
This becomes a dynamical system under the shift map
\[  \sigma\left((x_n)_{n\in \mathbb{Z}}\right) = (x_{n+1})_{n\in \mathbb{Z}}. \]
We introduce two subsystems of $\left(\left(\mathbb{R}/2\mathbb{Z}\right)^{\mathbb{Z}}, \sigma\right)$ by
\begin{equation*}
   \begin{split}
    \mathcal{Y} & := \left\{(x_n)_{n\in \mathbb{Z}}\in \left(\mathbb{R}/2\mathbb{Z}\right)^{\mathbb{Z}} \middle|\, 
    \forall n\in \mathbb{Z} : \max\left(\rho(x_n, x_{n+1}), \rho(x_{n+1}, x_{n+2})\right) = 1\right\}, \\
    \mathcal{Z} & := \left\{(x_n)_{n\in \mathbb{Z}}\in \left(\mathbb{R}/2\mathbb{Z}\right)^{\mathbb{Z}} \middle|\, 
    \forall n\in \mathbb{Z} : \max\left(\rho(x_n, x_{n+1}), \rho(x_{n+1}, x_{n+2})\right) \geq  \frac{1}{2} \right\}.
   \end{split}
\end{equation*}
We have $\mathcal{Y}\subset \mathcal{Z}$.
They are both fixed-point free dynamical systems.

The dynamical system $(\mathcal{Y},\sigma)$ is related to the marker property by the next lemma.
(The dynamical system $(\mathcal{Z},\sigma)$ will be used in the next section.)

\begin{lemma} \label{lemma: marker property and Y}
Let $(X, T)$ be a dynamical system having the marker property.
Then there is an equivariant continuous map from $X$ to $\mathcal{Y}$.
\end{lemma} 

\begin{proof}
We use Lindenstrauss' lemma (Lemma \ref{lemma: Lindenstrauss}) with $N=1$.
Then there is a continuous function $\varphi:X\to \mathbb{R}$ such that the set 
\[  E := \{x\in X|\, \varphi(Tx) \neq \varphi(x)+1\} \]
satisfies $E\cap T^{-1}E=\emptyset$.

Consider the composition of $\varphi$ and the natural projection $\mathbb{R}\to \mathbb{R}/2\mathbb{R}$.
We also denote it as $\varphi:X\to \mathbb{R}/2\mathbb{Z}$.

\begin{claim}  \label{claim: condition for Y}
 For any $x\in X$
\[  \max\left(\rho\left(\varphi(x), \varphi(Tx)\right), 
     \rho\left(\varphi(Tx), \varphi(T^2 x)\right)\right) = 1.  \]
\end{claim}

Indeed, if $x\not \in E$ then $\varphi(Tx) = \varphi(x)+1$ and hence 
$\rho\left(\varphi(x),\varphi(Tx)\right) =1$.
If $x\in E$ then $Tx\not \in E$. So $\varphi(T^2 x) = \varphi(Tx)+1$ and 
$\rho\left(\varphi(Tx),\varphi(T^2 x)\right) = 1$.
This proves the claim.

We define an equivariant continuous map $f:X\to \left(\mathbb{R}/2\mathbb{R}\right)^{\mathbb{Z}}$
by 
\[ f(x) := \left(\varphi(T^n x)\right)_{n\in \mathbb{Z}}. \]
Take any $x\in X$ and $n\in \mathbb{Z}$.
Applying the above Claim \ref{claim: condition for Y} to the point $T^n x$, we get 
\[ \max\left(\rho\left(\varphi(T^n x), \varphi(T^{n+1}x)\right), 
    \rho\left(\varphi(T^{n+1} x), \varphi(T^{n+2} x)\right)\right) = 1.  \]
Therefore $f(x)\in \mathcal{Y}$. So $f$ is an equivariant continuous map from $X$ to $\mathcal{Y}$.
\end{proof}

\section{Existence of a free dynamical system which does not have the marker property: 
Proof of Theorem \ref{theorem: marker property}} 
 \label{section: proof of theorem marker property}

In this section we combine all the preparations and construct a free dynamical system 
which does not have the marker property.

\subsection{From an infinite product to a finite product}  
\label{subsection: from infinite product to finite product}

Here we introduce a \textit{trick to reduce a problem on an infinite product to 
one on a finite product}.
 
Recall that we have introduced two dynamical systems $\mathcal{Y}$ and $\mathcal{Z}$
in \S \ref{subsection: dynamical systems Y and Z}.

\begin{lemma} \label{lemma: from infinite to finite}
Let $(X_m, T_m)$ $(m\in \mathbb{N})$ be a countable number of dynamical systems.
Consider their product\footnote{This is also a dynamical system. The map $\prod_{m=1}^\infty T_m$ is defined by 
\[  \left(\prod_{m=1}^\infty T_m\right) (x_1, x_2, x_3, \dots) = (T_1 x_1, T_2 x_2, T_3 x_3, \dots). \]}:
\[  \left(\prod_{m=1}^\infty X_m, \prod_{m=1}^\infty T_m\right). \]
Suppose there is an equivariant continuous map 
\[ f: \prod_{m=1}^\infty X_m \to \mathcal{Y}. \]
Then there are a natural number $M$ and an equivariant continuous map 
\[  g:X_1\times X_2\times \dots \times X_M \to \mathcal{Z}. \]
\end{lemma}

\begin{proof}
Set 
\[ (\mathbb{X}, \mathbb{T}) := \left(\prod_{m=1}^\infty X_m, \prod_{m=1}^\infty T_m\right). \]
For each natural number $M$ we set 
\[ (\mathbb{X}_M, \mathbb{T}_M) := \left(\prod_{m=1}^M X_m, \prod_{m=1}^M T_m\right). \]
We denote as $\pi_M:\mathbb{X}\to \mathbb{X}_M$ the natural projection.

Let $\varphi:\mathbb{X} \to \mathbb{R}/2\mathbb{Z}$ be the 0-th component of
$f:\mathbb{X}\to \mathcal{Y}\subset \left(\mathbb{R}/2\mathbb{Z}\right)^\mathbb{Z}$.
We have 
\[ f(x) = \left(\varphi\left(\mathbb{T}^n x\right)\right)_{n\in \mathbb{Z}}. \]

Since $\mathbb{X}$ is compact, $\varphi$ is uniformly continuous.
Hence we can choose a natural number $M$ such that for any two points $x, y\in \mathbb{X}$
\begin{equation} \label{eq: pi_M and small error}
    \pi_M(x) = \pi_M(y) \Longrightarrow \rho\left(\varphi(x), \varphi(y)\right) < \frac{1}{4}. 
\end{equation}

Fix a point $p\in X_{M+1}\times X_{M+2}\times X_{M+3}\times \cdots$.
We define a continuous map $\psi:\mathbb{X}_M\to \mathbb{R}/2\mathbb{Z}$ by 
\[  \psi(x) := \varphi(x, p). \]
We also define an equivariant continuous map 
$g:\mathbb{X}_M\to \left(\mathbb{R}/2\mathbb{Z}\right)^{\mathbb{Z}}$ by 
\[ g(x) := \left(\psi\left(\mathbb{T}_M^n x\right)\right)_{n\in \mathbb{Z}} =
   \left(\varphi\left(\mathbb{T}_M^n x, p\right)\right)_{n\in \mathbb{Z}}.   \]
We need to show $g(x)\in \mathcal{Z}$.

For any $x\in \mathbb{X}_M$ and $n\in \mathbb{Z}$
\[ \pi_M\left(\mathbb{T}_M^n x, p\right) = \pi_M\left(\mathbb{T}^n\left(x, p\right)\right)  \quad (= \mathbb{T}_M^n x).   \]
By (\ref{eq: pi_M and small error})
\begin{equation} \label{eq: error is small}
  \rho\left(\psi\left(\mathbb{T}_M^n x\right), \varphi\left(\mathbb{T}^n(x,p)\right)\right)
   = \rho\left(\varphi\left(\mathbb{T}_M^n x, p\right), \varphi\left(\mathbb{T}^n(x,p)\right)\right) < \frac{1}{4}.
\end{equation}

\begin{claim}  \label{claim: condition for Z}
For any $x\in \mathbb{X}_M$
\[  \max\left\{\rho\left(\psi(x), \psi\left(\mathbb{T}_M x\right)\right), 
     \rho\left(\psi\left(\mathbb{T}_M x\right), \psi\left(\mathbb{T}_M^2 x\right)\right)\right\}
 \geq \frac{1}{2}. \]
\end{claim}

Indeed, since $f(x,p)\in \mathcal{Y}$, there is $k\in \{0,1\}$ satisfying 
\[ \rho\left(\varphi\left(\mathbb{T}^k (x, p)\right), \varphi\left(\mathbb{T}^{k+1}(x,p)\right)\right) = 1. \]
Using (\ref{eq: error is small}) for $n=k$ and $n=k+1$,
\[ \rho\left(\psi\left(\mathbb{T}_M^k x\right), \varphi\left(\mathbb{T}^k(x,p)\right)\right) < \frac{1}{4}, \quad 
   \rho\left(\psi\left(\mathbb{T}_M^{k+1} x\right), \varphi\left(\mathbb{T}^{k+1}(x,p)\right)\right) < \frac{1}{4}. \]
By the triangle inequality 
\begin{equation*}
  \begin{split}
   \rho \left(\psi\left(\mathbb{T}_M^k x\right), \psi\left(\mathbb{T}_M^{k+1} x\right)\right)& \geq 
    \rho\left(\varphi\left(\mathbb{T}^k (x, p)\right), \varphi\left(\mathbb{T}^{k+1}(x,p)\right)\right) \\
    &-\rho\left(\psi\left(\mathbb{T}_M^k x\right), \varphi\left(\mathbb{T}^k(x,p)\right)\right) 
    -  \rho\left(\psi\left(\mathbb{T}_M^{k+1} x\right), \varphi\left(\mathbb{T}^{k+1}(x,p)\right)\right) \\
  & > 1 -\frac{1}{4} -\frac{1}{4} = \frac{1}{2}.
  \end{split}
\end{equation*}
This proves the above claim.

Take any point $x\in \mathbb{X}_M$ and any integer $n$.
We apply Claim \ref{claim: condition for Z} to the point $\mathbb{T}_M^n x$ and get 
\[  \max\left\{\rho\left(\psi\left(\mathbb{T}_M^n x\right), \psi\left(\mathbb{T}_M^{n+1} x\right)\right), 
     \rho\left(\psi\left(\mathbb{T}_M^{n+1} x\right), \psi\left(\mathbb{T}_M^{n+2} x\right)\right)\right\} \geq \frac{1}{2}. \]
This shows $g(x)\in \mathcal{Z}$. So $g$ is an equivariant continuous map from $\mathbb{X}_M$ to $\mathcal{Z}$.
\end{proof}

\subsection{Proof of Theorem \ref{theorem: marker property}} \label{subsection: proof of Theorem marker property}

Now we are ready to prove Theorem \ref{theorem: marker property}.
We write the statement again.

\begin{theorem}[$=$ Theorem \ref{theorem: marker property}]
There exists a free dynamical system which does not have the marker property.
\end{theorem}

\begin{proof}
Recall that the dynamical system $\mathcal{Z}$ 
(introduced in \S \ref{subsection: dynamical systems Y and Z})
is fixed-point free. 
Appying Corollary \ref{cor: increasing coindex by one} to $\mathcal{Z}$, we can choose a natural number $N$ and 
a positive number $\delta$ such that for all prime numbers $p$
\begin{equation}  \label{eq: choice of N and delta}
  \coind_p P_p\left(\mathcal{X}(N,\delta)\right) \geq \coind_p P_p(\mathcal{Z}) + 1.
\end{equation}

Now we consider the dynamical systems $\mathcal{X}_m(N,\delta)$ $(m\in \mathbb{N})$ introduced in 
\S \ref{subsection: variants of X(N,delta)}.
(The parameters $N$ and $\delta$ have been fixed by the condition (\ref{eq: choice of N and delta}).)
For simplicity of the notation, we set $\mathcal{X}_m := \mathcal{X}_m(N,\delta)$.
We define a dynamical system $X$ as the product of $\mathcal{X}_m$:
\[  X := \prod_{m=1}^\infty \mathcal{X}_m. \]
We are going to prove that $X$ is free and does not have the marker property.

For each natural number $m$, the system $\mathcal{X}_m$ has no $m$-periodic point
(Lemma \ref{lemma: coind of X_m(N,delta)}).
So $X$ has no $m$-periodic point. Since $m$ is arbitrary, $X$ is free.

Suppose $X$ has the marker property.
From Lemma \ref{lemma: marker property and Y}, there is an equivariant continuous map 
$f:X\to \mathcal{Y}$.

Appying Lemma \ref{lemma: from infinite to finite} to the infinite product $X = \prod_{m=1}^\infty \mathcal{X}_m$,
there are a natural number $M$ and an equivariant continuous map 
\[ g: \mathcal{X}_1\times \mathcal{X}_2\times \dots \times \mathcal{X}_M\to \mathcal{Z}. \]
Take a prime number $p$ larger than $M$, and restrict the map $g$ to the set of $p$-periodic points:
\[ g: P_p(\mathcal{X}_1)\times P_p(\mathcal{X}_2)\times \dots \times P_p(\mathcal{X}_M) \to P_p(\mathcal{Z}). \]
This is a $\mathbb{Z}_p$-equivariant continuous map.
So by Proposition \ref{prop: basic properties of coindex} (1)
\begin{equation}  \label{eq: coindex of product of X_m and Z}
     \coind_p\left\{P_p(\mathcal{X}_1)\times P_p(\mathcal{X}_2)\times \dots \times P_p(\mathcal{X}_M)\right\}
      \leq \coind_p P_p(\mathcal{Z}). 
\end{equation}
By Proposition \ref{prop: basic properties of coindex} (2), the left-hand side is equal to 
\[ \min\left\{\coind_p P_p(\mathcal{X}_1), \coind_p P_p(\mathcal{X}_2), \dots, \coind_p P_p(\mathcal{X}_M)\right\}. \]
By Lemma \ref{lemma: coind of X_m(N,delta)}, for all $1\leq m\leq M$
\[ \coind_p P_p(\mathcal{X}_m) = \coind_p P_p\left(\mathcal{X}(N,\delta)\right). \]
Hence 
\[ \min\left\{\coind_p P_p(\mathcal{X}_1), \coind_p P_p(\mathcal{X}_2), \dots, \coind_p P_p(\mathcal{X}_M)\right\}
     = \coind_p P_p\left(\mathcal{X}(N,\delta)\right). \]
Therefore the inequality (\ref{eq: coindex of product of X_m and Z}) becomes
\[  \coind_p P_p\left(\mathcal{X}(N,\delta)\right) \leq \coind_p P_p(\mathcal{Z}). \]
However, from (\ref{eq: choice of N and delta})
\[ \coind_p P_p\left(\mathcal{X}(N,\delta)\right) \geq \coind_p P_p(\mathcal{Z}) + 1. \]
This is a contradiction. Thus $X$ does not have the marker property.
\end{proof}

\begin{remark}
By investigating the above arguments a bit more closely, 
we can actually prove the following stronger statement:
\textit{Let $N\geq 1$ and $\delta>0$, and set 
\[  (X, T) := \prod_{m=1}^\infty \left(\mathcal{X}_m(N,\delta),\mathrm{shift}\right). \]
This is a free dynamical system.
If $N$ is sufficiently large and $\delta$ is sufficiently small, then there is no open set $U\subset X$
satisfying
\[ U \cap T^{-1}U = \emptyset, \quad X= \bigcup_{n\in \mathbb{Z}} T^{-n} U. \]}
\end{remark}

\section{Open problems}  \label{section: open problems}

This paper is just a starting point for investigating applications of 
$G$-index to dynamical systems theory.
There are definitely many open problems and new phenomena to be explored.
Here we mention just a few questions directly related to the main results of the paper.

Theorem \ref{theorem: first main result} states that,  for any fixed-point free dynamical system 
$(X,T)$, the sequence 
\[  \ind_p P_p(X) \quad (p=2,3,5,7,11,13,17,\dots) \]
grows at most linearly in $p$.
However we do not know whether there is an example which actually has linear growth.
Namely 

\begin{problem}  \label{problem: linear growth}
Is there a fixed-point free dynamical system $(X, T)$ satisfying 
\[  \ind_p P_p(X) \geq C p \]
for some positive constant $C$ and all sufficiently large prime numbers $p$?
We can ask the same question for $\coind_p P_p(X)$ and $\conn P_p(X)$.
\end{problem}

The authors have spent a lot of time trying to solve this problem.
But we have not succeeded.
It seems better to study the following simpler question before attacking Problem \ref{problem: linear growth}.

\begin{problem} \label{problem: unbounded}
Is there a fixed-point free dynamical system $(X, T)$ such that the sequence 
\[  \ind_p P_p(X) \quad (p=2,3,5,7,11,13,17,\dots) \]
is unbounded? We can ask the same question for $\coind_p P_p(X)$ and $\conn P_p(X)$.
\end{problem}

The above two problems might look abstract.
But indeed they are concrete questions.
By the universality of $\mathcal{X}(N,\delta)$ in Theorem \ref{theorem: universality},
the questions (concerning $\mathbb{Z}_p$-index and coindex\footnote{The quesiton about $\conn P_p(X)$
does not reduce to the case of $X=\mathcal{X}(N,\delta)$
because the connectivity does not necessarily increase under morphisms.}) reduce to the case of $X = \mathcal{X}(N,\delta)$.
So they are more or less equivalent to the following concrete problem.

\begin{problem}
Let $N$ be a natural number and $\delta$ a positive number.
For prime numbers $p$, estimate the $\mathbb{Z}_p$-index (or $\mathbb{Z}_p$-coindex) of 
\begin{equation*} 
   P_p\left(\mathcal{X}(N,\delta)\right) = 
   \left\{\left(x_0, x_1, \dots, x_{p-1}\right)\in \left([0,1]^N\right)^{\mathbb{Z}_p}\middle|\, 
       \forall n\in \mathbb{Z}_p: |x_n-x_{n+1}|\geq \delta\right\}.
\end{equation*}
Notice that we are mainly interested in the asymptotic behavior as $p\to \infty$ (while $N$ and $\delta$ are fixed).
\end{problem}

We also propose some problems on the marker property.
Theorem \ref{theorem: marker property} shows that periodic-point freeness alone does not imply 
the marker property.
So we need an additional condition besides the freeness for guaranteeing the marker property.
The following conjecture seems plausible.

\begin{conjecture}  \label{conjecture: zero mdim and marker property}
If a free dynamical system has zero mean dimension then it has the marker property.
\end{conjecture}

Indeed this is equivalent to the following conjecture of Lindenstrauss:

\begin{conjecture}[Lindenstrauss]  \label{conjecture: Lindenstrauss}
If a free dynamical system has zero mean dimension then it has the small boundary property.
\end{conjecture}

The \textit{small boundary property} is a dynamical analogue of totally disconnectedness 
introduced by Lindenstrauss--Weiss \cite[Section 5]{Lindenstrauss--Weiss}.
(Here we do not provide the precise definition.) 
Lindenstrauss \cite[Theorem 6.2]{Lindenstrauss} proved that 
if a \textit{zero mean dimensional} dynamical system has 
the marker property then it has the small boundary property.
On the other hand, Gutman (in a private communication) proved that if a free dynamical system has the small 
boundary property then it has the marker property.
So Conjectures \ref{conjecture: zero mdim and marker property} and 
\ref{conjecture: Lindenstrauss} are equivalent.

We notice that the dynamical system $\prod_{m=1}^\infty \mathcal{X}(N,\delta)$ constructed in 
\S \ref{subsection: proof of Theorem marker property} has infinite mean dimension.
It seems difficult to construct a finite mean dimensional example by our current method.
So we propose:

\begin{problem}
Construct a \textit{finite mean dimensional} free dynamical system which does not have the marker property.
\end{problem}

Finally we would like to mention a curious observation.
Let $\mathcal{Z}$ be the dynamical system introduced in \S \ref{subsection: dynamical systems Y and Z}.
\textit{If Conjecture \ref{conjecture: zero mdim and marker property} is true, then the sequence 
\[  \coind_p P_p(\mathcal{Z}), \quad (p=2,3,5,7,11,\dots) \]
must be unbounded.} 
(In particular, the answer to Problem \ref{problem: unbounded} becomes \textit{yes}
if Conjecture \ref{conjecture: zero mdim and marker property} is true.)

On the contrary, suppose there is a natural number $K$ satisfying
$\coind_p P_p(\mathcal{Z}) \leq K$ for all prime numbers $p$.

For each natural number $m$, we define a symbolic subshift $\Sigma_m\subset \{1,2,3\}^\mathbb{Z}$ by 
\[ \Sigma_m = \left\{(x_n)_{n\in \mathbb{Z}}  \in \{1,2,3\}^\mathbb{Z} \, \middle|\, 
                  \forall n\in \mathbb{Z}: x_n\neq x_{n+m}\right\}. \]
It is easy to check that $\Sigma_m$ has no $m$-periodic point and $P_p(\Sigma_m) \neq \emptyset$
for all prime numbers $p> m$.
We define a dynamical system $X_m$ by 
\[ X_m = \left(\Sigma_m\right)^{*(K+2)} \quad (\text{the join of the $(K+2)$ copies of $\Sigma_m$}). \]
The system $X_m$ has no $m$-periodic point and for prime numbers $p> m$
\begin{equation} \label{eq: coind_p of X_m is larger than K}
   \coind_p P_p(X_m) = \coind_p P_p(\Sigma_m)^{*(K+2)} \geq K+1 > \coind_p P_p(\mathcal{Z}). 
\end{equation}

Now we consider the infinite product $\prod_{m=1}^\infty X_m$. This is a free dynamical system.
By applying the argument of \S \ref{subsection: proof of Theorem marker property}, we can prove that 
$\prod_{m=1}^\infty X_m$ does not have the marker property. 
Indeed, suppose $\prod_{m=1}^\infty X_m$ has the marker property.
By Lemma \ref{lemma: marker property and Y} there is an equivariant continuous map 
$f:\prod_{m=1}^\infty X_m\to \mathcal{Y}$.
By Lemma \ref{lemma: from infinite to finite}, for some $M\geq 1$, there is an equivariant continuous map 
$g:\prod_{m=1}^M X_m \to \mathcal{Z}$.
Take a prime number $p>M$ and restrict $g$ to the set of $p$-periodic points:
\[  g: \prod_{m=1}^M P_p(X_m) \to P_p(\mathcal{Z}). \]
Then we get 
\[ \min_{1\leq m\leq M} \coind_p P_p(X_m) \leq \coind_p P_p(\mathcal{Z}). \]
This contradicts the above (\ref{eq: coind_p of X_m is larger than K}).

On the other hand, each $X_m$ is finite dimensional (and hence has zero mean dimension).
Then, by \cite[p. 231, (2.2)]{Lindenstrauss}, the product $\prod_{m=1}^\infty X_m$ has zero mean dimension.
If Conjecture \ref{conjecture: zero mdim and marker property} is true, then $\prod_{m=1}^\infty X_m$ has the marker property.
This is a contradiction.

Conjecture \ref{conjecture: zero mdim and marker property} (equivalently, 
Conjecture \ref{conjecture: Lindenstrauss}) is very difficult for our current technology.
However, the above system $\prod_{m=1}^\infty X_m$ has a very special form 
(i.e. infinite product of finite dimensional systems).
So we might be able to prove by our current technology that it has the marker property.

\vspace{0.5cm}

\address{ Masaki Tsukamoto \endgraf
Department of Mathematics, Kyushu University, Moto-oka 744, Nishi-ku, Fukuoka 819-0395, Japan}

\textit{E-mail}: \texttt{tsukamoto@math.kyushu-u.ac.jp}

\vspace{0.5cm}

\address{ Mitsunobu Tsutaya \endgraf
Department of Mathematics, Kyushu University, Moto-oka 744, Nishi-ku, Fukuoka 819-0395, Japan}

\textit{E-mail}: \texttt{tsutaya@math.kyushu-u.ac.jp}

\vspace{0.5cm}

\address{ Masahiko Yoshinaga \endgraf
Department of Mathematics, Hokkaido University, North 10, West 8, Kita-ku,
Sapporo, 060-0810, Japan}

\textit{E-mail}: \texttt{yoshinaga@math.sci.hokudai.ac.jp}


\begin{thebibliography}{99}





\bibitem[Gro99]{Gromov}
M. Gromov, 
Topological invariants of dynamical systems and spaces of holomorphic maps: I,
Math. Phys. Anal. Geom. \textbf{2} (1999) 323-415.




\bibitem[Gut15]{Gutman Jaworski}
Y. Gutman, Mean dimension and Jaworski-type theorems, 
Proceedings of the London Mathematical Society \textbf{111(4)} (2015) 831-850.


\bibitem[Gut17]{Gutman periodic points}
Y.~Gutman,
Embedding topological dynamical systems with periodic points in cubical shifts,
Ergodic Theory Dynam. System \textbf{37} (2017) 512-538.





\bibitem[GLT16]{Gutman--Lindenstrauss--Tsukamoto}
Y. Gutman, E. Lindenstrauss, M. Tsukamoto, 
Mean dimension of $\mathbb{Z}^k$-actions, 
Geom. Funct. Anal.
\textbf{26} Issue 3 (2016) 778-817.




\bibitem[GQT19]{Gutman--Qiao--Tsukamoto}
Y.~ Gutman, Y.~Qiao, M. Tsukamoto,
Application of signal analysis to the embedding problem of $\mathbb{Z}^k$-actions,
Geom. Funct. Anal.  \textbf{29} (2019) 1440-1502.




\bibitem[GT20]{Gutman--Tsukamoto}
Y.~Gutman, M.~Tsukamoto,
Embedding minimal dynamical systems into Hilbert cubes,
Invent. math.  \textbf{221} (2020) 113-166.



\bibitem[Lin99]{Lindenstrauss}
E. Lindenstrauss,
Mean dimension, small entropy factors and an embedding theorem,
Inst. Hautes \'{E}tudes Sci. Publ. Math. \textbf{89} (1999) 227-262.



 
\bibitem[LT19]{Lindenstrauss--Tsukamoto double VP}
E.~Lindenstrauss, M.~Tsukamoto,
Double variational principle for mean dimension, 
Geom. Funct. Anal., \textbf{29} (2019) 1048-1109.
 
 


\bibitem[LW00]{Lindenstrauss--Weiss}
E. Lindenstrauss, B. Weiss,
Mean topological dimension,
Israel J. Math. \textbf{115} (2000) 1-24.



 



\bibitem[Mat08]{Matousek}
J.~Matou$\mathrm{\check{s}}$ek,
Using the Borsuk--Ulam theorem,
2nd, corrected printing, Springer, 2008.




\bibitem[Tsu20]{Tsukamoto mdim potential}
M.~Tsukamoto,
Double variational principle for mean dimension with potential, 
Adv. Math. \textbf{361} (2020) 106935.


\end{thebibliography}
\end{document}